\documentclass[leqno,12pt]{article}
\usepackage{amssymb,amsmath,amsthm}
\def\eps{\varepsilon}
\def\osc{\hbox{osc}}
\def\f{\frac}
\def\st{such that}
\def\fy{\varphi}

\def\rest{\upharpoonright}

\def\colon{{:}\;}

\newfont{\bb}{msbm8 scaled\magstephalf}     

\newfont{\bbb}{msbm10 scaled\magstephalf}     
\def\N{\mbox{\bbb N}}

\def\R{\mbox{\bbb R}}

\def\lur{{\bf LUR}}

\def\trle{\preccurlyeq}
\def\trl{\prec}
\def\lsc{lower semicontinuous}

\def\1{{\mathchoice {\rm 1\mskip-4mu l} {\rm 1\mskip-4mu l}
{\rm 1\mskip-4.5mu l} {\rm 1\mskip-5mu l}}}

\def\arre{\par \medskip}
\def\half{{\textstyle\frac12}}

\newfont{\sssB}{msam10 scaled\magstephalf}     

\newtheorem{teorema}{Theorem}

\newtheorem{proposic}[teorema]{Proposition}
\newtheorem{lema}{Lemma}

\title{
\Huge {\bf Spaces of functions with countably many
discontinuities}}
\author{R. Haydon, A. Molt\'o\thanks{The second author has been
supported by BFM2003--07540/MATE, Ministerio de Ciencia y
Tecnolog\'{\i}a MCIT y FEDER (Spain).}, J. Orihuela\thanks{The third author has been
supported by Fundaci\'on S\'eneca CARM: 00690/PI/04 and by BFM2002--01719 del MCIT.}}
\date{}
\begin{document}

\maketitle

\begin{abstract}
Let $\Gamma$ be a Polish space and let $K$ be a separable and
pointwise compact set of functions on $\Gamma$.  Assume further that
each function in $K$ has only countably many discontinuities.  It is
proved that $\mathcal C(K)$ admits a $\frak T_p$-lower
semicontinuous and locally uniformly rotund norm, equivalent to the
supremum norm. A slightly more general result is shown and a related
conjecture is stated.
\end{abstract}

\begin{section}{Introduction}
It is known that, for a wide range of classes of ``well-behaved''
compact spaces $K$, the space $\mathcal C(K)$ of all continuous
real-valued functions on $K$ admits a norm, equivalent to the
supremum norm, that is locally uniformly rotund.  One approach to
this sort of problem, using projectional resolutions of the identity,
reached its most general result with the proof in \cite{Devilleetal} that
such a renorming exists when $K$ is a continuous image of a
Valdivia compact.  One of the present authors \cite{HaydonLUR} recently established LUR
renormability for $\mathcal C(K)$ when $K$ is a Namioka--Phelps
compact.  Novel techniques developed in \cite{transfer} led to a proof of
LUR renormability for the space $\mathcal C(H)$, where $H$ is the
Helly space.  The main result of this paper can be seen as a
generalization of this last result.  We use techniques from both
\cite{transfer} and \cite{HaydonLUR}, as well as borrowing a key idea from the
recent paper \cite{Argyrosetal} of Argyros et al.

The Helly space is an example of a class of  compacta about which
there are still a number of open questions.  We say that a compact
space $K$ is a Rosenthal compact if there exists a Polish space
$\Gamma$ and a homeomorphism from $K$ onto a subset of the space
$\mathcal B_1(\Gamma)$ of Baire--1 functions on $\Gamma$ equipped
with the pointwise topology.  Important (and now classical) results
of Rosenthal \cite{Rosenthal} and Bourgain, Fremlin and Talagrand
\cite{BFT} show that such compact spaces have certain properties
which place them close to metrizable compacta and to weakly compact
subsets of Banach spaces. It is natural therefore to hope for good
results when we look at these compact spaces in the context of
renorming theory.  In fact, however, as Todorcevic \cite{Tod1} has
recently observed, there is a scattered Rosenthal compactification
$K$ of a tree space such that C(K) has no LUR renorming,
\cite{HaydonTrees}. Now that space $K$ is non-separable and, as
other recent work of Todorcevic \cite{Tod2} has shown, it is only
from separable Rosenthal compacta that we should expect really good
behaviour.  We are therefore led to make the conjecture.

\medskip\noindent
{\bf Conjecture.}
If $K$ is a separable Rosenthal compact then $\mathcal C(K)$
admits a locally uniformly convex renorming.
\medskip

A proof of this conjecture would yield as an immediate corollary
that $X^*$ is LUR renormable whenever $X$ is a separable Banach
space with no subspace isomorphic to $\ell^1$.  Indeed, in this
case, we may take $\Gamma$ to be the unit ball of the dual space
$X^*$, which is compact and  metrizable (so certainly Polish)
under the weak* topology $\sigma(X^*,X)$ and
$K$ to be the unit ball of $X^{**}$ under the weak* topology
$\sigma(X^{**},X^*)$.
By the results of \cite{OdellRosenthal}, the elements of $K$ are
then of the first Baire class when we regard them as functions on
$\Gamma$.  Moreover, $K$ is separable, since the unit ball of $X$
(which we are assuming to be separable) is dense in $K$ by
Goldstine's theorem. Finally, of course, $X^*$ embeds as a closed
subspace of $\mathcal C(K)$.
The main theorem of this paper establishes LUR renormability of $\mathcal C(K)$
only for a subclass of separable Rosenthal compacta $K$, namely those
representable as spaces of functions with only countably many
discontinuities in the Polish space $\Gamma$.

\begin{teorema}
[Main Theorem]\label{MainThm}
Let $\Gamma$ be a Polish space and let $K$ be a separable and pointwise compact
set of functions on $\Gamma$.  Assume further that each function
in $K$ has only countably many discontinuities.  Then $\mathcal C(K)$
admits a $\frak T_p$-lower semicontinuous and  locally uniformly
rotund norm, equivalent to the supremum norm.
\end{teorema}
\end{section}

When $K$ is not assumed to be separable and each element of $K$ has
only countably many discontinuities the fact that $C(K)$ is a
$\sigma$--fragmentable space for the pointwise convergence topology
${\frak T}_p$ was obtained by I. Kortezov \cite{Kortezov} following
previous results on the Namioka property by A. Bouziad
\cite{Bouziad}. In the last section of the paper we present our
approach for this result. We believe it could help to deal with the
LUR renorming problem on this class of $C(K)$ spaces.

This paper had its origin in a visit of the second and third authors
to the University of  Oxford in 2002, they wish to express their
gratitude to  Brasenose College, to the Mathematical Institute and
to the members of the Functional Analysis Group for their kind
hospitality.

\begin{section}
{Preliminaries}

Our notation and terminology are standard: we write $\omega$ for
the set $\{0,1,2,\dots\}$ of all natural numbers and $\N$ for the
set of all positive integers.  When $A$ is a set, we write $\# A$
for the cardinality of $A$ and $[A]^{<\omega}$ for the set of all
finite subsets of $A$. We recall that a topological space is said
to be Polish if it separable, metrizable and complete for some
metric compatible with the topology.  A space is analytic if it is
a continuous image of some Polish space.

We shall be considering real-valued functions on a Polish space
$\Gamma$ which have only countably many discontinuities. For such
a function $s$ we may introduce the following
subsets of $\Gamma$, which consist of the ``big'' discontinuities:
\begin{align*}
J(s,\delta) &= \{\gamma\in \Gamma \colon \osc \ s\restriction
U>\delta\text{\ whenever\ }U\text{\ is open and\ }\gamma\in U\},
\end{align*}
the above being defined for all positive real $\delta$.
Each of these sets is a countable closed subset of the Polish
space $\Gamma$ and hence a scattered topological space.  We recall the
Cantor-Bendixon derivation for such spaces:  as usual, for a
topological space $J$, we write $J'$ for the derived set,
consisting of those points of $J$ which are not isolated; by
transfinite recursion, we define successive derived sets
$J^{(\xi)}$ for ordinals $\xi$ by
$$
J^{(0)}=J, \qquad\qquad
J^{(\eta)}=\bigcap_{\xi<\eta}(J^{(\xi)})'\quad (\eta>0).
$$
The space $J$ is scattered if and only if
$J^{(\lambda)}=\emptyset$ for some $\lambda$; if this is so, then
the smallest such $\lambda$ is called the derived length, or
Cantor-Bendixon index of $J$.  Of course, if $J$ is countable,
this index will automatically be a countable ordinal.

 The proof of our main theorem
divides into two parts. First we establish LUR-renormability of
$\mathcal C(K)$ without making a separability assumption about
$K$, assuming instead that there is some countable ordinal
$\Omega$ such that $J(s,\delta)^{(\Omega)}=\emptyset$ for all
$\delta>0$ and all $s\in K$.  Then we show that separability of
$K$ implies such a uniform bound on derived length.  It is in the
second of these steps that we use ideas from \cite{Argyrosetal},
involving the rank of a well-founded relation, and the so-called
Rank Theorem for analytic relations.

In our case, the relation will be the ordering on a certain tree. We
now set out some notation and terminology.  Recall that a tree is a
partially ordered set $(\Upsilon, \trl)$ with the property that for
each $v\in \Upsilon$ the set $\{u\in \Upsilon \colon u\trl v\}$ is
well-ordered by $\trl$.  Any tree $\Upsilon$ may be partitioned into
{\bf levels} $\Upsilon_\xi$, with $\Upsilon_0$ consisting of the
$\trl$-minimal elements.  More generally, for $v\in \Upsilon$ we
define the {\bf height} $\text{ht}(v)$ to be the order-type of
$\{u\in \Upsilon \colon u\trl v\}$, and define the $\xi^{\text{th}}$
level $\Upsilon_\xi$ to be $\{v\in \Upsilon \colon
\text{ht}(v)=\xi\}$. Following Todorcevic \cite{TodTrees}, when
$\Upsilon$ and $\Psi$ are trees, we write $\Upsilon\otimes\Psi$ for
the set $\{(u,v)\in\Upsilon\times\Psi \colon
\text{ht}(u)=\text{ht}(v)\}$, which is itself a tree when equipped
with the order $(u,v)\trle (u',v')$ if and only if $u\trle u'$ and
$v\trle v'$. (The restriction that $\text{ht}(u)=\text{ht}(v)$ is
important to ensure that we have a total order on the predecessors
of $(u',v')$.)

The trees with which we shall be concerned are {\bf well-founded},
that is to say that they contain no strictly ($\trl$)-increasing
sequence.  Of course this implies that it is only the levels
$\Upsilon_n$ ($n\in \omega$) which are non-empty, but it also
enables to introduce a second ordinal index, which we call the {\bf
rank}.  An element of a well-founded tree is of rank 0 if it is
{\bf$\trl$-maximal}.  More generally, we define a derivation on
$\Upsilon$, by setting $A' =A\setminus\text{max\,}A$, where
$\text{max\,}A$ is the set of all maximal elements of the subset $A$
of $\Upsilon$; we then define
$$
\Upsilon^{[0]}=\Upsilon;\qquad
\Upsilon^{[\eta]}=\bigcap_{\xi<\eta} {\Upsilon^{[\xi]}}'
\quad(\eta\ge 1).
$$
An element $u$ of $\Upsilon$ is defined to be of rank $\xi$ if
$u\in \Upsilon^{[\xi]}\setminus \Upsilon^{[\xi+1]}$.  It is
straightforward to check that the rank $r(u)$ satisfies the
following identity
$$
r(u) = \sup\{r(v)+1 \colon v\in u^+\},
$$
where $u^+$ is the set of immediate successors of $u$ in the
tree-order $\trl$ (and where we are of course taking the supremum
of the empty set of ordinals to be 0).
It will be useful to record here an easy result about ranks in the
tree $\Upsilon\otimes\Upsilon$.

\begin{lema}\label{otimes}
Let $\Upsilon$ and $\Psi$ be a well-founded tree.  Then
$\Upsilon\otimes\Psi$ is well-founded and the rank of $(u,v)$
in $\Upsilon\otimes\Psi$ is given by
$$
r(u,v) = \min\{r(u),r(v)\}.
$$
\end{lema}
\begin{proof}
The immediate successors of $(u,v)$ in $\Upsilon\otimes\Psi$
are the pairs $(u',v')$ with $u'\in u^+$ and $v'\in v^+$.  So the
result follows from the identity we noted above.
\end{proof}

The important theorem that we shall use is the Rank Theorem.  We
refer the reader to \cite{Kechris} for a full account of and state
here the special case which concerns us.

\medskip\noindent
{\bf The Rank Theorem.} {\it
Let $\Upsilon$ be a well-founded tree and assume that the tree
order $\trle$ (considered as a subset of $\Upsilon\times
\Upsilon$)
is an analytic topological space.  Then there is a countable ordinal
$\Omega$ such that $\Upsilon^{[\Omega]}=\emptyset.$
}

Finally we recall the definition of locally uniform convexity and
some (probably familiar) convexity arguments.
If $\nu$ is a non-negative real-valued convex function on a real
vector space $X$ and $x, x_n\in X$ $(n\in \N)$, we say that the
LUR hypothesis holds for $\nu$ (and $ x$ and the sequence $x_n$)
if $\nu(x_n)$ and $\nu(\half(x+x_n))$ both
tend to the limit $\nu(x)$ as $n\to \infty$.  If $\nu$ is a norm and
the LUR hypothesis implies that $\nu(x-x_n)\to 0$, then we say that
$\nu$ is locally uniformly rotund.

It is worth noting that the LUR hypothesis holds
if and only if
$$
\half \nu(x)^2+\half\nu(x_n)^2-\nu(\half(x+x_n))^2\to 0
$$
as $n\to \infty$.
If $\nu^2=\mu^2+\lambda^2$, where $\mu$ and $\lambda$ are
themselves non-negative and convex, and if the LUR hypothesis
holds for $\nu$, then it holds for $\lambda$ and $\mu$ also.  This
is an observation that we shall use a number of times, justifying
each such application by the phrase ``by convexity''.

\end{section}

\begin{section}
{Construction of a LUR norm}

In this section we shall prove the following theorem.

\begin{teorema}
Let $\Gamma$ be a Polish space and let $K$ be a pointwise compact
set of Baire-1 functions on $\Gamma$.  Assume that there exists a countable ordinal $\Omega$
such that for all $s\in K$ and all $\delta>0$  the
$\Omega^{\text{th}}$ derived set $J(s,\delta)^{(\Omega)}$ is empty. Then the Banach space $\mathcal C(K)$ admits an
equivalent $\mathfrak T_p$-lower semicontinuous and locally uniformly rotund norm.
\label{LURthm}
\end{teorema}

As in \cite{HaydonTrees} and \cite{HaydonLUR} we shall employ a
method of recursive definitions, combined with the following
result, which we refer to as Deville's Lemma.    Let us mention that an approach based on countably decomposition in the spirit of \cite{transfer} is also possible.   Notice that, using the language
introduced in the previous section,  the key assumption in this lemma
can be expressed by saying that the LUR hypothesis holds for
$\theta$, $x$ and $(x_n)$.

\begin{lema}[\cite{DGZ}, p.279]
Let $(\phi_i)_{i\in I}$ and $(\psi_i)_{i\in I}$ be two pointwise-bounded
families of
non-negative, real-valued, convex functions on a real vector space
$Z$.    For $i\in I$ and positive integers $p$ define
functions $\theta_{i,p}$, $\theta_p$ and $\theta$ by setting
\begin{align*}
2\theta_{i,p}(x)^2 &= \phi_i(x)^2 + p^{-1}\psi_i(x)^2\\
\theta_p(x) &= \sup_{i\in I}\theta_{i,p}(x)\\
\theta(x)^2 &= \sum_{p=1}^\infty 2^{-p}\theta_p(x)^2.
\end{align*}
Let $x$ and $x_r$ $(r\in \N)$ be  elements of $Z$
and assume that
$$
\half \theta(x)^2+\half\theta(x_r)^2-\theta(\half(x+x_r))^2\to 0
$$
as $r\to \infty$.  Then there is a sequence $(i_r)$ of elements of
$I$ such that
\begin{align*}
\phi_{i_r}(x)&\to \sup_{i\in I} \phi_i(x)\quad\text{and}\\
\half\psi_{i_r}(x)^2 &+
\half\psi_{i_r}(x_r)^2-\psi_{i_r}(\half(x+x_r))^2\to 0
\end{align*}
as $r\to \infty$.
\end{lema}

Modifying a little our earlier notation for discontinuity sets,
when $s,t\in K$ and $m\in \mathbb N$, we shall write $J(s,t,m)$
for the union $J(s,1/4m)\cup J(t,1/4m)$.  Our hypothesis implies
that $J(s,t,m)^{(\Omega)}=\emptyset$ for all $s,t\in K$.
Let $\mathcal B$ be a countable base for the topology in $\Gamma$.
Let $Q$ be a dense countable subset of $\Gamma$, and write
$P=\Gamma \setminus Q$.  When $R\subset Q$, $F\subset P$ and
$m\in \N$,  we define
$$
I(R,F,m ) = \left\{ (s,t)\in K \times K  \colon
\|(s-t)\rest_{R}\|_{\infty} \leq 1/4m,\ \|(s-t)\rest_{F}\|_{\infty}
\leq 1/m \right\} .
$$
When $R,F,m$ are as above, $\xi <\Omega$ is an ordinal and
$\mathcal M$ is a finite subset of $\mathcal B$, we define
$$
I(R,F,m,\xi ,\mathcal M )=\{(s,t)\in I(R,F,m ) \colon
\#J(s,t,m)^{(\xi )}\cap U \le 1\text{\ for\ all\ }U\in\mathcal M\}.
$$
These sets $I(R,F,m,\xi,\mathcal M)$
will play the role of the index set $I$ in applications of
Deville's Lemma.  In order to make clear our applications of this lemma
we make the (otherwise redundant) definition, which defines the
functions $\phi_i(x)=\fy (x,s,t)$ for $i=(s,t)\in I$,
$$
\fy
(x,s,t)=(1/2)|x(s)-x(t)|, \ s,t\in K,\ x\in C(K),
$$
and introduce the suprema
$$
\fy (x,R,F,m )=\sup \{ \fy (x,s,t) \colon  (s,t)\in I(R,F,m ) \}
$$
and
$$
\fy (x,R,F,m,\xi ,\mathcal M )=\sup \{ \fy (x,s,t) \colon  (s,t)\in
I(R,F,m,\xi ,\mathcal M) \} .
$$

The definition of the convex functions $\psi_i$ is more complicated
and some more notation is needed.    Without loss of generality we
can assume that $K\subset [0,1]^{\Gamma}$. For finite $R\subset Q$
and $F\subset P$ we choose a finite subset $K(R,F,m)$ of $K$ (of
cardinality at most $m^{\# F}(4m)^{\# R}$) \st\ for all $s\in K$
there exists $t\in K(R,F,m)$ with $(s,t)\in I(R,F,m)$.  To construct
the functions $\psi_i$ and the promised equivalent norm, we make
recursive definitions as set out in the following lemma.\newpage

\begin{lema}
There are functions $\nu,\theta, \psi$, defined for $x\in C(K)$,
$R\subset Q$, $F\subset P$, $m\in \N$, $\xi <\Omega$, $\mathcal M\in
[\mathcal B]^{<\omega}$, $k \in \N$, $s$, $t\in I(R,F,m,\mathcal M
)$, and satisfying the following:

\begin{align*}
3\nu (x,R,F,m )^{2} =& \fy (x,R,F,m )^{2}+
\f{1}{\# K(R,F,m)}\sum_{t\in K(R,F,m)}x(t)^{2}\\
 &\quad +\sum_{\xi ,\mathcal M}a( \xi ,\mathcal
M)\theta ( x,R,F,m,\xi ,\mathcal M )^{2}\\
\theta ( x,R,F,m,\xi ,\mathcal M )^{2} =&
\sum_{k=1}^{\infty}2^{-k}\ \theta(x,R,F,m,\xi ,\mathcal M ,k )^{2};\\
2 \theta(x,R,F,m,\xi ,\mathcal M ,k )^{2} =& \sup_{(s,t)\in I(Q,F,m,\xi
,\mathcal M )} \left[ \fy(x,s,t)^{2}+k^{-1}
\psi (x,s,t,R,F,m,\xi ,\mathcal M )^{2}\right] ; \\
 \psi (x,s,t,R,F,m,\xi ,\mathcal M )^{2} =& \f{1}{\# \mathcal
M}\sum_{U\in \mathcal M}\nu (x,R,F\cup (U\cap J(s,t,m)^{(\xi)}) ,m )^{2} .
\end{align*}
where the positive constants $a(\xi ,\mathcal M )$ are chosen so
that $$\sum_{\xi<\Omega ,\ \mathcal M\in [\mathcal B]^{<\omega}}a(\xi ,\mathcal M )=1.$$
\end{lema}
\begin{proof}
As in \cite{HaydonTrees} and \cite{HaydonLUR}, this follows from
Banach's Contraction Mapping Theorem applied in a suitable function
space.  To describe it we take as domain for our functions the set
$\mathcal{D} : =C(K)\times [Q]^{<\omega}\times [R]^{<\omega}\times
\N \times [0,\Omega )\times [\mathcal{B}]^{<\omega}\times \N$ and we
consider the set $\mathcal{H}$ of all mappings from $\mathcal{D}$
into $\R^2$, i.e. $(\nu ,\theta ) (x,R,F,m,\xi ,\mathcal{M},k )
=\left( \nu (x,R,F,m ),\theta (x,R,F,m,\xi ,\mathcal{M},k ) \right)$
\st\ $\nu (\cdot ,R,F,m )$ and $\theta (\cdot ,R,F,m,\xi
,\mathcal{M},k )$ are convex and positively homogeneous in $x\in
C(K)$ and they both are bounded above by $\|\cdot \|_{\infty}$.  On
$\mathcal{H}$ we consider the complete metric given by $d\left(
\left( \nu, \theta\right) ,\left( \nu^{\prime} ,\theta^{\prime}
\right)\right) =\max \left\{\sup \left\{ \left| \nu^{2}(\delta )-
\nu^{\prime 2}(\delta )\right| \colon  \delta \in
\mathcal{D}\right\} ,\sup \left\{ \left| \theta^{2}(\delta )-
\theta^{\prime 2}(\delta )\right|  \colon  \delta \in
\mathcal{D}\right\} \right\}$ and we define the contractive map
$F\colon \mathcal{H}\to \mathcal{H}$, $F(\nu, \theta )=(\nu^{\prime}
, \theta^{\prime}  )$ where
\begin{eqnarray*}
3\nu^{\prime} (x,R,F,m )^{2} &=&\fy (x,R,F,m )^{2}+
\f{1}{\# K(R,F,m)}\sum_{t\in K(R,F,m)}x(t)^{2}+\\[12pt]
 &&\quad +\sum_{\xi ,\mathcal M}a( \xi ,\mathcal M)
 \sum_{k=1}^{\infty}2^{-k}\ \theta(x,R,F,m,\xi ,\mathcal M ,k )^{2};\\[12pt]
2 \theta^{\prime}(x,R,F,m,\xi ,\mathcal M ,k )^{2} &=&
\sup_{(s,t)\in I(Q,F,m,\xi ,\mathcal M )} \bigg[ \fy(x,s,t)^{2}+\\[12pt] &&\qquad +k^{-1}
 \f{1}{\# \mathcal M}
\sum_{U\in \mathcal M}\nu (x,R,F\cup (U\cap J(s,t,m)^{(\xi)}) ,m )^{2} \bigg]  .
\end{eqnarray*}
\end{proof}

Notice that,
because of the definition of $I(Q,F,m,\xi
,\mathcal M )$, each of the sets $F\cup (U\cap J(s,t,m)^{(\xi)})$
which occur in the definition of $\psi$ is either just $F$, or else $F$
with one extra element appended. Notice also that $\nu$, $\theta$ and $\psi$, when considered
as functions of $x$, are non-negative, positively homogeneous and
bounded (by 1) on the unit ball of $\mathcal C(K)$.

We define a pointwise \lsc\ norm in $\mathcal C(K)$
by
$$ 2\| x\|^{2}=\| x\|^{2}_{\infty}+\sum_{m,R}  c(R,m)\nu (x,
R,\emptyset ,m )^{2}$$ where $c(R,m)$ are further positive constants
\st\ $\sum_{m,R} c(R,m)=1$. Our aim is to show that this
norm is \lur .

So we consider $x\in \mathcal C(K)$ and a sequence $x_n$ such that
the LUR hypothesis is satisfied for this norm. Notice that, by
convexity and the definition of the norm as an $\ell^2$ sum, the
LUR hypothesis holds for each of the functions
$\nu(\cdot,R,\emptyset,m)$.

Now let $\eps >0$ be given.  It will be enough to show that there
exists a subsequence $\left( x_{n_{k}}\right)$ \st\ $\left\|
x-x_{n_{k}}\right\|_{\infty} <5\eps$ for all $k$.

Our compact space $K$ is a closed subset of $[0,1]^\Gamma$
equipped with the product topology and our given function $x$ is continuous on $K$.
Hence, there exist $m\in \N$ and a finite subset $T$ of $\Gamma$ \st\ $|x(s)-x(t)| \le \eps$
whenever $\sup_{\gamma\in T}|s(\gamma)-t(\gamma)|\le 1/m$.  If we set $S=T\cap P$
we obviously have $|x(s)-x(t)| \le\eps$ whenever $(s,t)\in I(Q,S,m)$.
Rather than working with this set $S$, we choose a finite
subset $S$ of $P$ of minimal cardinality subject to the above
condition, that is
$$
(s,t)\in I(Q,S,m)\Longrightarrow |x(s)-x(t)| \le\eps .
$$
Now, by an easy compactness argument, we may also choose a finite
subset $R$ of $Q$ \st\ $|x(s)-x(t)| <2\eps $
whenever $ (s,t)\in I(R,S,m)$.  Recalling the definitions given
earlier, we see that
$$
\fy(x,R,S,m) <\eps .
$$
The proof of Theorem~\ref{LURthm}
depends on two lemmas.  It is perhaps worth emphasizing that from
now on $x$, $x_n$, $\eps$, $m$, $R$ and $S$ are all fixed
\begin{lema}
If the LUR hypothesis holds for the function $\nu(\cdot,R,S,m)$
and a subsequence $(x_{n_k})$ then
$\left\| x-x_{n_k}\right\|_{\infty} <5\eps $ for all large
enough $k$. \label{A}
\end{lema}
\begin{proof}
By convexity and the expression for $\nu$ as an $\ell^2$-sum, we get
$$
\f{1}{2}\fy (x,R,S,m)^{2} +\f{1}{2}\fy \left( x_{n},R,S,m
\right)^{2}-\fy \left( \f{1}{2}\left(x+x_{n}\right) ,R,S,m
\right)^{2}\rightarrow 0 \ \hbox{as}\ n\to \infty.
$$
Since $\fy (x,R,S,m) <\eps$, it follows that
$$
\fy \left( x_{n},R,S,m \right) <\eps \ \hbox{for all large enough}\
n.
$$
Looking at the second term in the definition of $\nu$ and applying convexity again we see that
$$
\f{1}{2}x(t)^{2}
 +\f{1}{2}x_{n}(t)^{2}
-\left( \f{1}{2}\left(x(t)+x_{n}(t)\right) \right)^{2}\rightarrow
0 \ \hbox{as}\ n\to \infty \ \hbox{for any}\ t\in
K(R,S,m),
$$
which implies that
$$
\max_{t\in
K(R,S,m)}\left| x(t)-x_{n}(t)\right| <\eps \text{\ for\ all\ large\
enough\ }n.
$$

For any $s\in K$ there exists $t\in K(R,S,m)$ with
$(s,t)\in I(R,S,m)$ and so
$$
\left| x(s)-x_{n}(s)\right|  \leq
2\fy ( x,R,S,m) + 2\fy \left( x_{n},R,S,m\right) +\max_{t\in
K(R,S,m)}\left| x(t)-x_{n}(t)\right| <5\eps
$$
for all large enough $n$.
\end{proof}

\begin{lema}
Let $F$ be a proper subset of $S$, let $(x_{n_k})$ be a subsequence of $(x_n)$
and assume that the LUR hypothesis holds for the subsequence $(x_{n_k})$
and the function $\nu(\cdot,R,F,m)$.
Then there exists $\gamma\in S\setminus F$ and a further subsequence
$(x_{n_{k_j}})$ such that the LUR hypothesis holds for $(x_{n_{k_j}})$
and the function $\nu(\cdot,R,F\cup\{\gamma\},m)$.
\label{B}
\end{lema}
\begin{proof}
 To simplify notation, avoiding double subscripts,
let us assume that the initial subsequence $(x_{n_k})$ is actually the sequence $(x_n)$ itself.
Since $\# F <\# S$ we may use the minimality in the
choice of $S$ to see that there is some
$$
(s,t)\in I(Q,F,m) \quad \hbox{with}\quad  |x(s)-x(t)| > \eps .
$$
By the choice of $S$, there exists $\gamma\in S$ with $|s(\gamma
)-t(\gamma )| >(1/m)$ . We claim that any such $\gamma$ must be in
$J(s,t,m)$, so that
$S \cap  J(s,t,m) \neq \emptyset .$
Indeed, otherwise, $\gamma\notin J(s,t,m)$ and there exists an open
set $U\ni \gamma$ \st\ $| s(\delta )-s(\gamma) | \le1/4m$ and $|
t(\delta )-t(\gamma) | \le1/4m$ for all $\delta \in U$.  By density
of $Q$, there is some $\delta \in Q\cap U$ and, since $(s,t)\in
I(Q,F,m)$ we have $|s(\delta )-t(\delta )| \le1/4m$. Combining the
above inequalities, we obtain $|s(\gamma )-t(\gamma )| \le 3/4m$, a
contradiction. \arre

Since $S$ is finite there is some maximum
ordinal $\xi(s,t,m)$ with
$$
S\cap  J(s,t,m)^{(\xi (s,t,m))}
\neq \emptyset .
$$
We now assume that we have chosen
$(s,t)\in  I(R,F,m)$ in such a way that
$$
\xi(s,t,m) =\min \left\{ \xi \left( s^{\prime},t^{\prime},m \right)
 \colon  \left( s^{\prime},t^{\prime} \right) \in I(Q,F,m) \text{ and }
|x(s^{\prime})-x(t^{\prime})|>\eps\right\}.
$$
Let $\mathcal M$ be a finite subset of $\mathcal B$, chosen in such a
way that $S\subseteq \bigcup \mathcal M$ and $\#U\cap J(s,t,m)^{(\xi )}
\le 1$ for all $U\in \mathcal M$.

We have $(s,t)\in I(Q,F,m,\xi ,\mathcal M)$ and $\fy (x,s,t)>\eps/2$, so
$$
\fy (x,Q,F,m,\xi ,\mathcal M ) > \eps/2 .
$$
If we now look at the third term in the definition of $\nu$ and apply
the familiar convexity argument we see that
the function $\theta ( \cdot ,R,F,m,\xi ,\mathcal M)$ and the sequence
$(x_{n})$ satisfy the LUR hypothesis.
 So by Deville's
Lemma there is a sequence $\left( s_{n},t_{n} \right) \in
I(Q,F,m,\xi ,\mathcal M )$ \st\
$$
\f{1}{2}\psi \left( x,s_{n},t_{n},R,F,m \right) ^{2} +\f{1}{2}\psi
\left( x_{n},s_{n},t_{n},R,F,m \right)^{2}- \psi \left(
\f{1}{2}\left(x+x_{n}\right) ,s_{n},t_{n},R,F,m
\right)^{2}\rightarrow 0
$$
$$\text and\quad
\fy \left( x,s_{n},t_{n} \right) \rightarrow \fy (x,Q,F,m,\xi
,\mathcal M )>\eps/2 .
$$
So $\fy \left( x,s_{n},t_{n} \right) >\epsilon/2$,
i.e.$ |x(s_n)-x(t_n)|>\epsilon$  for all large enough $n$.
Reasoning as before,
we get that $ S \cap J\left(s_{n},t_{n},m\right) \neq \emptyset $
for such $n$.  From the minimality of $\xi$ it follows
that $ S \cap J\left(s_{n},t_{n},m\right)^{(\xi )} \neq
\emptyset $. Now
$$
 S \cap J\left(s_{n},t_{n},m\right)^{(\xi )}\subseteq
 \bigcup_{U\in \mathcal M}U\cap J\left(s_{n},t_{n},m\right)^{(\xi )}
 $$
and each of the sets $U\cap J\left(s_{n},t_{n},m\right)^{(\xi )}$
contains at most one element because $(s_n,t_n)\in I(Q,F,m,\xi ,\mathcal M
)$.  Thus, for some $U(n)\in \mathcal M$ the intersection
$U(n)\cap J\left(s_{n},t_{n},m\right)^{(\xi )}$ contains exactly one point
$\gamma_{n}$ which is in $S$. If we proceed to a
subsequence $(n_k)$, we may assume that $U_{n_k}$ and $\gamma_{n_{k}}$ are the same
set $U$ in ${\mathcal M}$ and the same element $\gamma$ of $S$, for all $k$. Finally, by looking at the
definition of $\psi $ and applying convexity yet again, we see that
$$
\f{1}{2}\nu (x,R,F\cup \{ \gamma \},m)^{2} +\f{1}{2}\nu \left(
x_{n_{k}},R,F\cup \{ \gamma \},m \right)^{2}-\nu \left(
\f{1}{2}\left(x+x_{n_{k}}\right) ,R,F\cup \{ \gamma \},m
\right)^{2}
$$
\begin{eqnarray*}
&&  \leq \#\mathcal M \left( \f{1}{2}\psi (x,s_{n_{k}},t_{n_{k}},R,F,m)
^{2}+ \f{1}{2}\psi \left( x_{n_{k}},s_{n_{k}},t_{n_{k}},R,F,m
\right)^{2} \right. \\[12pt]
&&~\hspace{4cm}\left. -\psi \left( \f{1}{2}\left(x+x_{n_{k}}\right)
,s_{n_{k}},t_{n_{k}},R,F,m \right)^{2}\right) \rightarrow 0 \
\hbox{as}\ k\rightarrow \infty .
\end{eqnarray*}
This is the LUR hypothesis for $\nu(\cdot,R,F\cup\{\gamma\},m)$
and the subsequence $(x_{n_k})$.
\end{proof}

To complete the proof, we have already noted that the LUR hypothesis
holds for $(x_n)$ and the function $\nu(\cdot,R,\emptyset,m)$.
A finite number of applications of Lemma~\ref{B} yield a subsequence satisfying
the LUR hypothesis for the function $\nu(\cdot, R,S,m)$.  Lemma~\ref{A} now
does what we need to complete the proof of Theorem~\ref{LURthm}.
\end{section}

\begin{section}{Boundedness of the Cantor--Bendixon index}

We devote this section to a proof of the following result, which,
together with Theorem~\ref{LURthm} clearly yields Theorem~\ref{MainThm}.  Our
proof is closely modeled on a recent theorem of Argyros et al
\cite{Argyrosetal}.

\begin{teorema}\label{index}
Let $\Gamma$ be a Polish space and let $K$ be a pointwise compact
set of functions on $\Gamma$.  Assume that each function
$s\in K$ has only countably many discontinuities and that the set
$K$ is separable.  Then there exists a countable ordinal $\Omega$
such that for all $s\in K$ and all $\delta>0$  the
$\Omega^{\text{th}}$ derived set $J(s,\delta)^{(\Omega)}$ is empty.
\end{teorema}

We first establish a little notation: when $k$ is a natural number
we write $D_k$ for the set of all finite sequences of 0's and 1's,
of length at most $k$.   Thus $D_k = \bigcup_{0\le j\le k} \{
0,1\}^j$.  The set $D_0$ has just one element, the sequence $()$ of
length 0.  If $\sigma \in \{0,1\}^j$ the {\bf length} $j$ of
$\sigma$ will be denoted $|\sigma|$ and we shall write $\sigma.0$
(resp. $\sigma.1$) for the element of $\{0,1\}^{j+1}$ which extends
$\sigma$ and has $0$ (resp. 1) in its last place.  We define
$0.\sigma$ and $1.\sigma$ analogously.

We fix a metric $d$ on $\Gamma$, compatible with the given
topology, under which $\Gamma $ is complete, as well as a
countable base $\mathcal B$ for the topology of $\Gamma$ and a sequence
$(s_m)$ which is dense in $K$.  We now introduce, for
all infinite subsets $M$ of $\mathbb N$, all natural numbers
$p$, all nonempty open subsets $X$ of $\Gamma$, and all
positive real numbers $\delta$, a set $\Upsilon(X,M,p,\delta)$.
Elements of this set are tuples
$$
(M,k,(U_\sigma)_{\sigma\in D_k},(\alpha_\sigma)_{\sigma\in
D_k},(\beta_\sigma)_{\sigma\in D_k} ),
$$
where $k\in \omega$,
$U_\sigma\in \mathcal B$, $\alpha_\sigma\in \Gamma$ and $\beta_\sigma\in \Gamma$  satisfy the
following conditions:
\begin{enumerate}
\item $\overline {U_\sigma}\subseteq X$ and $\text{diam\,}U_\sigma\le 2^{-p-|\sigma|}$ for all $\sigma\in D_k$;
\item $\overline {U_{\sigma.0}}\cap\overline {U_{\sigma.1}}=\emptyset$ and
$\overline {U_{\sigma.0}}\cup\overline {U_{\sigma.1}}\subseteq U_\sigma$
whenever $\sigma \in D_{k-1}$;
\item $\alpha_\sigma\,,\, \beta_{\sigma}\in U_\sigma$;
\item $|s_m(\alpha_{\sigma})-s_m(\beta_{\sigma})|>\delta$ for all large enough $m\in M$.
\end{enumerate}
The union $\bigcup_{M\in [\mathbb N]^{\omega}}\Upsilon(X,M,p,\delta)$
will be denoted $\Upsilon(X,p,\delta)$.  This is a tree under the
following ordering:
$$
(M,k,(U_\sigma)_{\sigma\in D_k},(\alpha_\sigma)_{\sigma\in D_k} ,(\beta_\sigma)_{\sigma\in D_k} )
\prec (M',k',(U'_\sigma)_{\sigma\in
D_{k'}},(\alpha'_\sigma)_{\sigma\in D_{k'}},(\beta'_\sigma)_{\sigma\in D_{k'}})
$$
if and only if
$$
M=M'\,,\quad k < k'\,,\quad U_\sigma=U'_\sigma \
,\ \alpha_\sigma = \alpha'_\sigma
\ \text{and}\ \beta_\sigma =
\beta'_\sigma\ \quad\text{for\ all\ }\sigma\in D_k.
$$
Notice that the coordinate $k$ equals
the height of the element
$(M,k,(U_\sigma) ,(\alpha_\sigma) ,(\beta_\sigma))$ in the tree.

\begin{lema}
For each $\delta>0$ and each non-empty open subset $X$ of $\Gamma$,
the tree $\Upsilon(X,p,\delta)$ is well founded.
\end{lema}
\begin{proof}
The reader will probably have realized that elements of the tree
$\Upsilon(X,p,\delta)$ can be regarded as finite attempts at
constructing a Cantor set of discontinuities for some element of
$K$.  The proof of the present lemma makes this idea more
explicit.

We have to show that our tree has no infinite branch.  So suppose,
if possible, such a branch exists.  It would consist of a sequence of elements
$$
(M,k,(U_\sigma) ,(\alpha_\sigma),(\beta_\sigma) )\quad (k\ge 1)
$$
 satisfying (1) to  (4).

It follows from completeness of $\Gamma$ and the conditions (1) and (2) that,
for each infinite sequence $z\in
\{0,1\}^\omega$, the intersection $\bigcap_{l\in \omega}
U_{z\restriction l}$ contains just one point $\gamma_z$.  If $s$ is any element of $K$ which is a
cluster point of $(s_m)_{m\in M}$ then we have
$$
|s(\beta_\sigma)-s(\alpha_\sigma)|\ge \delta
$$
for all $\sigma\in \{0,1\}^{<\omega}$.
Since each of the sequences $(\alpha_{z\restriction l})$ and $(\beta_{z\restriction
l})$
converges to $\gamma_z$, we see that each $\gamma_z$ is a discontinuity
point of $s$, contrary to our hypothesis that there are only countably many such discontinuities.
\end{proof}

\begin{lema}
For each non-empty open subset $X$ of $\Gamma$, each $p\in \omega$ and each
$\delta>0$, the relation $\trl$ on $\Upsilon(X,p,\delta)$ is analytic.
\end{lema}
\begin{proof}
Here we follow \cite{Argyrosetal} quite closely.  We note that $X$
is open in the Polish space $\Gamma$ and hence itself Polish, and
that $[\omega]^\omega$ is also a Polish space. The countable sets
$\omega$ and $\mathcal B$ will be equipped with the discrete
topology.  So if we define
$$
H= \bigcup_{k\in \omega}
\left([\omega]^{\omega}\times\{k\}\times \mathcal B^{D_k}
\times X^{D_{k}}\times X^{D_{k}}\right),
$$
$H$ is a disjoint union of Polish spaces and hence a Polish space.
It follows from our description of the elements of
$\Upsilon(X,p,\delta)$ that $\Upsilon(X,p,\delta)\subseteq H$.
We shall show that the relation $\trl$ is an analytic subset of
$H\times H$.  Now it is very easy to see that $\trl$ is closed in
$\Upsilon(X,p,\delta)\times\Upsilon(X,p,\delta)$, so it will be enough
for us to show that $\Upsilon(X,p,\delta)$ is analytic.

Now it is a standard result that we may enhance the
topology of the Polish space $\Gamma$ in such a way that all the sets $U\in
\mathcal B$ are both open and closed, and all the functions $s_m$
are continuous, while $\Gamma$ remains a Polish space.  Let us write
$\tilde X$ for $X$ equipped with this enhanced topology, and
$\tilde H$ for $H$ equipped with a similarly enhanced topology.
What we shall show is that $\Upsilon(X,p,\delta)$ is a countable union of closed
subsets of $\tilde H$.

We set
\begin{align*}
\Upsilon_n&=\{(M,k,(U_\sigma) ,(\alpha_\sigma)_{\sigma\in
D_k},(\beta_\sigma))_{\sigma\in
D_k} )\in\Upsilon(X,p,\delta) \colon \\
&\qquad|s_m(\alpha_\sigma)-s_m(\beta_\sigma)|\ge\delta \text{
whenever } \sigma\in D_k\text{ and }n\le m\in M\}
\end{align*}
and shall show that each $\Upsilon_n$ is closed in $\tilde H$.
Suppose then that
$$
(M^l,k^l,(U^l_\sigma)_{\sigma\in
D_{k^l}},(\alpha^l_\sigma)_{\sigma\in D_{k^l}},(\beta^l_\sigma)_{\sigma\in
D_{k^l}})\quad(l\in \omega)
$$
is a sequence of elements of $\Upsilon_n$ which converges
in $\tilde H$ to the limit
$$
(M,k,(U_\sigma)_{\sigma\in D_{k}} ,(\alpha_\sigma)_{\sigma\in D_{k}} ,(\beta_\sigma)_{\sigma\in
D_{k}}).
$$
Since we have equipped $\omega$ and $\mathcal B$ with
the discrete topology,  $k^l=k$ and $U^l_\sigma=U_\sigma$ for all
large enough $l$.  Since $\alpha_\sigma^l$ and $\beta_\sigma^l$
converge respectively to $\alpha_\sigma$ and $\beta_\sigma$ in the
topology of $\tilde X$, and since $U_\sigma$ is closed in that
topology, we have $\alpha_\sigma,\beta_\sigma\in U_\sigma$.  If
$m\in M$ and $m\ge n$ then $m\in M^l$ for all large enough $l$ so
that  $|s_m(\alpha^l_\sigma)-s_m(\beta^l_\sigma)|\ge \delta$ for all
large enough $l$.  Since we have arranged for $s_m$ to be
continuous in the topology of $\tilde X$, we see that
$|s_m(\alpha_\sigma)-s_m(\beta_\sigma)|\ge \delta$.  We have
finished showing that $\Upsilon_n$ is  closed
in $\tilde H$.
\end{proof}

To finish the proof of our theorem, we need to show that the rank
of the tree $\Upsilon(X,p,\delta)$ dominates the derived length of the set
$J(s,\delta)$ when $s\in K$.  We do this using two final lemmata,
the first of which expresses an obvious idea in what is perhaps
over-pedantic notation.

\begin{lema}\label{XY}
Let $X$ be a non-empty open subset of $\Gamma$, let $M$ be an infinite
subset of $\omega$, let $p$ be a
natural number and let $\delta$ be a positive real number.  Let
$Y$ and $Z$ be disjoint non-empty open subsets of $X$ and assume that there
exists $U\in \mathcal B$ with $Y\cup Z\subseteq U\subseteq X$ and
$\text{{\rm diam}}\,U\le 2^{-p}$.  If $\Upsilon(Y,M,p+1,\delta)^{[\xi]}$ and
$\Upsilon(Z,M,p+1,\delta)^{[\xi]}$ are both non-empty, then
$\Upsilon(X,M,p,\delta)^{[\xi+1]}$ is also non-empty.
\end{lema}
\begin{proof}
We shall show how to embed the tree
$\Upsilon(Y,M,p+1,\delta)\otimes \Upsilon(Z,M,p+1,\delta)$ into the
set of non-minimal elements of $\Upsilon(X,M,p,\delta)$, that is
to say, the elements of height at least 1. Our hypothesis,
together with Lemma~\ref{otimes} will then tell us that $\Upsilon(X,M,p,\delta)^{[\xi]}$
contains a non-minimal element, which in turn implies that
$\Upsilon(X,M,p,\delta)^{[\xi+1]}\ne\emptyset$.

Since $Y\subseteq U$ and $\Upsilon(Y,M,p+1,\delta)\ne \emptyset$
there certainly exist $\alpha, \beta\in U$ with
$|s_m(\alpha)-s_m(\beta)|\ge \delta$ for all large enough $m\in M$.
We define  $U_{()}=U$, $\alpha_{()}=\alpha$ and $\beta_{()}=\beta$.
Now let $(M,k,(V_\sigma)_{\sigma\in D_k},(\kappa_\sigma)_{\sigma\in
D_k},(\lambda_\sigma)_{\sigma\in D_k})$ and $(M,k,(W_\sigma)_{\sigma\in D_k},(\mu_\sigma)_{\sigma\in
D_k},(\nu_\sigma)_{\sigma\in D_k})$ be height k elements of
$\Upsilon(Y,M,p+1,\delta)$ and $ \Upsilon(Z,M,p+1,\delta)$
respectively.  If we define
\begin{align*}
U_{0.\sigma}=V_\sigma \qquad&\qquad U_{1.\sigma}=W_\sigma\\
\alpha_{0.\sigma}=\kappa_\sigma \qquad&\qquad \alpha_{1.\sigma}=\mu_\sigma\\
\beta_{0.\sigma}=\lambda_\sigma \qquad&\qquad \beta_{1.\sigma}=\nu_\sigma,
\end{align*}
then it is easy to check that $(M,k+1,(U_\sigma)_{\sigma\in D_{k+1}},
(\alpha_\sigma)_{\sigma\in D_{k+1}},(\beta_\sigma)_{\sigma\in
D_{k+1}})$ is a height $k+1$ element of $\Upsilon(X,M,p,\delta)$.
This defines the promised embedding.
\end{proof}

\begin{lema}
Let $X$ be a non-empty open subset of $\Gamma$ and let $\delta$ be a positive real number.  Let
$s\in K$ be the pointwise limit $s=\lim_{M\ni m\to \infty}s_m$
along the subsequence $M$, and assume that $X\cap
J(s,\delta)^{(\xi)}\ne \emptyset.$  Then, for every $p\in \omega$,
$\Upsilon(X,M,p,\delta)^{[\xi]}\ne \emptyset.$
\end{lema}
\begin{proof}
We proceed by induction on the ordinal $\xi$, starting with
$\xi=0$: if $J(s,\delta)\cap X \ne \emptyset$ we choose any
$\gamma$ in this set and then select $U\in \mathcal B$ with
$\gamma\in U$, $\overline U\subseteq X$, $\text{diam}\,U\le
2^{-p}$.  Since $\gamma\in J(s,\delta)$ we have
$\text{osc}(s\restriction U)>\delta$ so we can choose $\alpha,
\beta\in U$ with $|s(\alpha)-s(\beta)|>\delta.$  If we set
$U_{()}=U,\ \alpha_{()}=\alpha,\ \beta_{()}=\beta,$ then
$(M,0,(U_\sigma)_{\sigma\in D_0},(\alpha_\sigma)_{\sigma\in D_0},(\beta_\sigma)_{\sigma\in
D_0})\in \Upsilon(X,M,p,\delta)$.

Now suppose that $X\cap
J(s,\delta)^{(\xi+1)}\ne \emptyset $ and that the result is true
for $\xi$.  Choose $\gamma\in X\cap
J(s,\delta)^{(\xi+1)}$ and an element $U$ of $\mathcal B$ with
$\gamma\in U$, $\text{diam}\,U\le 2^{-p}$.  Since $\gamma$ is a
limit point of $J(s,\delta)^{(\xi)}$, we may find $\zeta\in U\cap
J(s,\delta)^{(\xi)}$ with $\zeta\ne \gamma$, and
then choose disjoint open $Y$ and $Z$ containing $\gamma$ and
$\zeta$ respectively.  By our inductive hypothesis (which of
course applies to {\it all} $p$), $\Upsilon(Y,M,p+1,\delta)^{[\xi]}$ and
$\Upsilon(Z,M,p+1,\delta)^{[\xi]}$ are both non-empty.  By
Lemma~\ref{XY} we now have $\Upsilon(X,M,p,\delta)^{[\xi+1]}\ne \emptyset.$

Finally, let $\eta$ be a limit ordinal, with
$X\cap J(s,\delta)^{(\eta )}\ne \emptyset $, and assume that the
result is true for all $\xi<\eta$.  As before, we choose
$\gamma\in X\cap J(s,\delta)^{(\eta)}$, $U\in \mathcal B$ with
$\gamma\in U$, $\text{diam}\,U\le 2^{-p}$ and $\alpha, \beta\in U$
with $|s_m(\alpha)-s_m(\beta)|\ge \delta$ for all large enough
$m\in M$.  We use these to define a height 0 element $(M,0,U,\alpha,\beta)$ of
$\Upsilon(X,M,p,\delta)$.  For any $\xi<\eta$, $\gamma$ is a limit
point of $X\cap J(s,\delta)^{(\xi)}$, so we can find disjoint open
subsets $Y,Z$ of $U$ such that $Y\cap J(s,\delta)^{(\xi)}$ and
$Z\cap J(s,\delta)^{(\xi)}$ are both nonempty. By inductive hypothesis,
both  $\Upsilon(Y,M,p+1,\delta)^{[\xi]}$ and
$\Upsilon(Z,M,p+1,\delta)^{[\xi]}$ are both non-empty. The proof of
Lemma~\ref{XY} shows that our already constructed height 0 element
$(M,0,U,\alpha,\beta)$ is in $\Upsilon(X,M,p,\delta)^{[\xi+1]}$.
Since this is true for all $\xi<\eta$ we have
$\Upsilon(X,M,p,\delta)^{[\eta]}\ne \emptyset$ as claimed.
\end{proof}

The proof of Theorem~\ref{index} is now complete.

\end{section}

\begin{section}{$\sigma$--Fragmentability of $C(K)$ for the non--separable case}

Let us recall that space $C(K)$ with the pointwise topology ${\frak
T}_p$ is said to be {\em $\sigma$--frag\-mentable} by its norm if
for every $\eps >0$ we can decompose
$C(K)=\bigcup_{n=1}^{\infty}C_{n,\eps }$ in such a way that for
every $n\in \N$ and every non--empty subset $T\subset C_{n,\eps }$
there exists a ${\frak T}_p$--open set $V$ \st\ $V\cap T$ is
non--empty and has norm diameter less than $\eps$.  This notion was
introduced and studied in \cite{JNRPL} where among other things it
is proved that $C(K)$ is $\sigma$--fragmentable if for every $\eps
>0$ we can decompose $C(K)=\bigcup_{n=1}^{\infty}C_{n,\eps }$ in
such a way that for every $n\in \N$ and every non--empty subset
$T\subset C_{n,\eps }$ there exists a ${\frak T}_p$--open $V$ \st\
$V\cap T$ is non--empty and covered by countably many sets of
diameter less than $\eps$, see \cite[Theorem~4.1.]{JNRPL}.

In our last section we shall prove the following
\begin{teorema}
Let $\Gamma$ be a Polish space and let $K$ be a
pointwise compact set of functions on $\Gamma$ \st\ each function $s\in K$ has only countably many discontinuities.  Then $(C(K),{\frak T}_{p})$ is $\sigma$--fragmentable by its norm.
\label{frag1824}
\end{teorema}
We first establish a little notation:  For $s\in K$ let $J(s)=\bigcup_{\delta >0}J(s,\delta )$, i.e.
the set of all the discontinuity points of $s$.  A finite sequence $\left\{
\left(s_{i},t_{i}\right)\right\}_{i=1}^{n}$,
$\left(s_{i},t_{i}\right)\in K\times K$, $1\leq i\leq n$, is said to
be {\em fitted} whenever for every $i$, $1\leq i\leq n$, we have
\begin{equation}
s_{i}(\gamma ) =t_{i}(\gamma )\ \hbox{for all}\ \gamma \in Q\cup
\left( \bigcup_{j<i}\left( J\left(s_{j}\right)\cup J\left(
t_{j}\right) \right)\right) , \label{frag1832}
\end{equation}
and the fitted sequence $\left\{
\left(s_{i},t_{i}\right)\right\}_{i=1}^{n}$ is said to have {\em
length} $n$. \arre Given $\eps >0$ we say that $x\in C(K)$
$\eps${\em --jumps} a fitted sequence $\left\{
\left(s_{i},t_{i}\right)\right\}_{i=1}^{n}$ whenever
\begin{equation}
\left| x\left( s_{i}\right) -x\left( t_{i}\right) \right| >\eps \ \hbox{for every}\ 1\leq i\leq n. \label{frag1833}
\end{equation}
Theorem~\ref{frag1824} will follow from some lemmata.
\begin{lema}
Given $x\in C(K)$ and $\eps >0$ there
exists $n\in \N$ \st\ $x$ does not $\eps$--jump any fitted sequence
in $K\times K$ which has length strictly bigger than $n$.
\label{frag1943}
\end{lema}
\begin{proof} Since $K$ is compact $x$ is uniformly continuous and there
exists $\delta >0$ and a finite subset $F$ of $\Gamma$ \st\
\begin{equation}
\left( \sup_{\gamma \in F}|s(\gamma )-t(\gamma )| < \delta
\Longrightarrow |x(s)-x(t)| <\eps \right) \ \hbox{for all}\ s,t\in
K. \label{frag1831}
\end{equation} Let $n : =\# F$, we claim that if $\left\{
\left(s_{i},t_{i}\right)\right\}_{i=1}^{m}$ is a fitted sequence
which is $\eps$--jumped by $x$ then we have $m\leq n$.  Indeed since
$Q$ is dense in $\Gamma$ from (\ref{frag1832}) it follows that
$s_{1}(\gamma )=t_{1}(\gamma )$ for every point of continuity
$\gamma$ of $s_1$ and $t_1$.  Then from (\ref{frag1832}),
(\ref{frag1833}) and (\ref{frag1831}) it follows
$$
F\cap \left( J\left( s_{1}\right)\cup J\left( t_{1}\right)\right)
\setminus Q\neq \emptyset .
$$
An obvious induction argument gives that
\begin{equation}
F\cap  \left( J\left(s_{i}\right)\cup J\left( t_{i}\right) \right)
\setminus \left( Q\cup \left( \bigcup_{j< i}\left(
J\left(s_{j}\right)\cup J\left( t_{j}\right) \right) \right) \right)
\neq \emptyset ,\qquad 1\leq i\leq m .\label{frag1850}
\end{equation}
Now the statement follows from (\ref{frag1850}). \end{proof} \arre It might be worth remarking that (\ref{frag1850}) shows that if a fitted sequence $\left\{
\left(s_{i},t_{i}\right)\right\}_{i=1}^{n}$ is $\eps$--jumped by some $x\in C(K)$ then $  \left(s_{i},t_{i}\right)\neq \left(s_{j},t_{j}\right)$ for
$1\leq i\neq j\leq n$.
\arre Given
$x\in C(K)$ and $\eps >0$ let $j(x,\eps )$ the minimum of the
natural numbers for which the thesis of Lemma~\ref{frag1943} holds
for $\eps$ and $x$.  Thus for any $j>j(x,\eps )$ the function $x$ cannot $\eps$--jump any fitted sequence of length $j$, and there exists a fitted sequence $\eps$--jumped by $x$ of length equal to $j(x,\eps )$ whenever $j(x,\eps ) >0$.

\arre Given $x\in C(K)$ and $\eps >0$, a subset
$S\subset \Gamma$ is said to $\eps$--{\em control} $x$ whenever
there exist a finite subset $F\subset S$ and $\delta >0$ for which
(\ref{frag1831}) holds.
\begin{lema}
Given $x\in C(K)$ and $\eps >0$ if $\left\{
\left(s_{i},t_{i} \right) \right\}_{i=1}^{j(x,\eps /2)}$ is a fitted
sequence which is $(\eps /2)$--jumped by $x$ then $Q\cup \left(
\bigcup_{i\leq j(x,\eps /2 )}\left( J\left(s_{i}\right)\cup J\left(
t_{i}\right)\right) \right)$ $\eps$--controls $x$. \label{frag1200}
\end{lema}
\begin{proof} Otherwise for any finite set $F$, $F\subset Q\cup \left(
\bigcup_{i\leq j(x,\eps /2 )}\left(  J\left(s_{i}\right)\cup J\left(
t_{i}\right) \right) \right)$, and any $n\in \N$ we can choose
$s(F,n)$, $t(F,n)\in K$ \st\ $$|s(F,n)(\gamma ) -t(F,n)(\gamma )|
<1/n\ \hbox{for all} \ \gamma \in F\ \hbox{whereas}\
|x(s(F,n))-x(t(F,n))|\geq \eps .$$  Since $K$ is compact there must
exist an adherent point $\left( \tilde{s},\tilde{t}\right)$ to the
net $\left( s(F,n),t(F,n)\right)$.  Now we have
$$
\left| x\left( \tilde{s}\right)-x\left( \tilde{t}\right)\right|\geq
\eps
>\eps /2 \ \hbox{and}\ \tilde{s}(\gamma )=
\tilde{t}(\gamma ) \ \hbox{for all}\ \gamma \in Q\cup \left(
\bigcup_{i\leq j(x,\eps /2 )} \left( J\left(s_{i}\right)\cup J\left(
t_{i}\right) \right) \right) .
$$
Then $\left(s_{1},t_{1} \right)$, \ldots, $\left(s_{j(x,\eps /2
)},t_{j(x,\eps /2 )} \right)$, $\left( \tilde{s}, \tilde{t}\right)$
is a fitted sequence of length $1+j(x,\eps /2)$ which $x$ $(\eps
/2)$--jumps, a contradiction. \end{proof}
\begin{lema}
For every $\eps >0$ there exists a
decomposition $C(K)=\bigcup_{n\in \N}C_{n,\eps }$ \st\ for every $n\in
\N$ and any $x\in C_{n,\eps }$ there exist a weak open set $W$ and a
countable set $N\subset \Gamma$ \st\ $x\in W\cap C_{n,\eps }$ and
$N$ $\eps$--controls every $y\in W\cap C_{n,\eps }$.
\label{frag1131}
\end{lema}
\begin{proof} Let $\eps >0$.  Let $C_{n,\eps }:=\left\{ x\in C(K) \colon  j(x,\eps
/2)=n\right\}$.  According to the proof of Lemma~\ref{frag1200} if
$n=0$ then the set $Q$ $\eps$--controls every function $y \in
C_{0,\eps }$.  Suppose $n>0$, and let us fix a fitted sequence
$\left\{ \left(s_{i},t_{i}\right)\right\}_{i=1}^{n}$ which is $\eps
/2$--jumped by $x$. Let $\delta >0$ \st\
$$
\left| x\left(s_{i}\right)-x\left(t_{i}\right)\right|
>\f{\eps}{2}+\delta ,\qquad 1\leq i\leq n.
$$
Then every element $y$ in the weak open set
$$
W:=\left\{ y \in C(K) \colon  \left|
x\left(s_{i}\right)-y\left(s_{i}\right)\right| <\delta /2 \
\hbox{\&}\ \left| x\left(t_{i}\right)-y\left(t_{i}\right)\right|
<\delta /2 ,\ 1\leq i\leq n\right\}
$$
$\eps /2$--jumps the above fitted sequence.  Then from
Lemma~\ref{frag1200} we conclude that the set
$$Q\cup \left( \bigcup_{i\leq j(x,\eps /2 )}\left(
J\left(s_{i}\right)\cup J\left( t_{i}\right)  \right) \right)$$
$\eps$--controls every $y\in C_{n,\eps }\cap W$. \end{proof}
\par The proof of Theorem~\ref{frag1824} now follows from the remarks in the beginning of this section together with the above lemma and the next proposition which is nothing else than a quantitative version of Ascoli's theorem.
\begin{proposic}
Let $K$ be a compact space embedded in a cube $[0,1]^{\Gamma}$ and let ${\mathcal F}$ be a bounded subset of $C(K)$.  If $\mathcal F$ is $\eps$--equicontinuous, i.e. there is a finite set $F\subset \Gamma$ and $\delta >0$ \st\ $|x(s)-x(t)| <\eps$ whenever $|s(\gamma )-t(\gamma )| <\delta$ for all $\gamma \in F$ and x in $\mathcal F$, then $\mathcal F$ can be covered by finitely many sets of norm diameter less than $3\eps$.
\label{frag1826}
\end{proposic}
\begin{proof}  Let us assume that $\mathcal F$ is bounded by $m$ and let us split up the interval $[-m,m]$ into a
finite number of sets of diameter less than $\eps $,
\begin{equation}
[-m,m] =\bigcup_{i=1}^{\ell}I_{i};\quad \hbox{diam}~\left(
I_{i}\right) <\eps , \ 1\leq i\leq \ell .
\label{framu820}
\end{equation}
Given a subset $S\subset \Gamma $ the symbol $\pi_{S} \colon
[0,1]^{\Gamma }\to [0,1]^S$ will stand for the canonical projection.
For the finite set $F$ we can split $\pi_{F}(K)$ up into finitely
many sets which are included into a cartesian product of intervals
of length strictly less than $\delta $, i.e.
$$
\pi_{F}(K) =\bigcup_{j=1}^{p}\left( L_{j}\cap \pi_{F}(K)\right)
$$
where $L_{j}\cap \pi_{F}(K)\neq \emptyset$ and each $L_j$ is a cube
whose factor intervals have a length strictly less than $\delta $.
Choose $s_{1}$,\ldots,$s_{p}\in K$ \st\ $\pi_{F}\left( s_{j}\right)
\in L_j$, $1\leq j\leq p$.  Then we can cover
$$
{\mathcal F} =  \bigcup_{h=1}^{N}T_{h}
$$
in such a way that each $T_{h}$ has the property that $y\left(
s_{j}\right)$ belongs to the same interval $I_{i}$ from (\ref{framu820}) for all $y\in T_{h}$
and all $j$, $1\leq j\leq p$.  \arre  We claim that the diameter of
every $T_{h}$ is less or equal than $3\eps$.   Indeed, given $y_1$, $y_{2}\in T_{h}$ and $s\in K$ there exists $j$, $1\leq
j\leq p$, \st\ $\pi_{F}(s)\in L_j$.  From the choice of $L_j$ we get
$$
\left| s(\gamma )-s_{j}(\gamma )\right| <\delta \ \hbox{for all}\
\gamma \in F.
$$
Now from our condition of $\eps$--equicontinuity it follows that
\begin{equation}
\left| y_{i}(s)-y_{i}\left( s_{j}\right)\right| < \eps , \
i=1,2. \label{frag1251}
\end{equation}
From the choice of $T_{h}$ we have that $y_{1}\left( s_{j}\right)$
and $y_{2}\left( s_{j}\right)$ belong to the same interval
$I_{i}$ from (\ref{framu820}) so
\begin{equation}
\left| y_{1}\left( s_{j}\right) -y_{2}\left( s_{j}\right)\right|
< \eps . \label{frag1253}
\end{equation}
From (\ref{frag1251}) and (\ref{frag1253}) we conclude that
$$
\left| y_{1}\left( s\right) -y_{2}\left( s\right)\right| <3\eps ,
$$
being the reasoning valid for every $s\in K$ the norm--diameter of
$T_{h}$ is not bigger than $3\eps$. \end{proof}

To complete the proof of Theorem \ref{frag1824} let us observe that
the decomposition from lemma \ref{frag1131} gives us sets $C_{n,\eps
}$ such that for every $x \in C_{n,\eps }$ the weak open set $W$
containing $x$ gives us the set $W\cap C_{n,\eps }$ which is
a countable union of $\eps$- equicontinuous sets, so it can be
covered by countably many sets of norm diameter less than $3\eps$ by
the former proposition. The conclusion now follows from
\cite[Theorem~4.1.]{JNRPL}

\end{section}

\arre \noindent {\sc R. Haydon: Brasenose College, Oxford OX1 4AJ,
England}\arre \noindent {\em email address: {\tt
richard.haydon@brasenose.oxford.ac.uk}} \arre \noindent {\sc A.
Molt\'o: Departamento de An\'alisis Matem\'atico, Facultad de
Matem\'a\-ticas, Universidad de Valencia, Dr. Moliner 50, 46100
Burjasot (Valencia), Spain}\arre \noindent {\em email address: {\tt
anibal.molto@uv.es}}\arre \noindent {\sc J. Orihuela: Departamento
de Matem\'aticas, Universidad de Murcia, Campus de Espinardo, 30100
Espinardo, Murcia, Spain} \arre \noindent {\em email address: {\tt
joseori@um.es}}
\end{document}